\documentclass{CSML}
\pdfoutput=1

\usepackage{lastpage}
\newcommand{\dOi}{13(3:16)2017}
\lmcsheading{}{1--\pageref{LastPage}}{}{}%
{May~10, 2016}{Aug.~25, 2017}{}

\newcommand{\orduparrow}{\mathord{\uparrow}}

\usepackage{enumerate}
\usepackage{hyperref}
\hypersetup{hidelinks}

\usepackage[]{amsmath, amssymb}

\begin{document}

\title[First Order Theories of Some Lattices of Open Sets]{First Order Theories of Some Lattices of Open Sets}

\author{Oleg Kudinov}
\address{S.L. Sobolev Institute of Mathematics
\\Siberian Branch of the
Russian Academy of Sciences\\
    Novosibirsk, Russia} 
\email{kud@math.nsc.ru}
 \author{Victor Selivanov}
\address{A.P. Ershov Institute of
Informatics Systems\\Siberian Branch of the Russian Academy of
Sciences, Novosibirsk, and Kazan (Volga Region) Federal University\\
    Russia}
    \email{vseliv@iis.nsk.su}

\keywords{Topological space, lattice, open set, effectively open set, first order theory, decidability, $m$-reducibility, interpretation.}
 \subjclass{2012 ACM CCS: [Mathematics of computing]: Continuous mathematics— Topology— Point-set topology; Continuous mathematics— Continuous functions; [Theory of computation]: Models of computation— Computability— Turing machines. 2010 Mathematics Subject Classification: 03D78, 03D45, 03D55, 03D30}

\begin{abstract}
We show that the first order
theory of the lattice of open sets in some natural topological spaces is
$m$-equivalent to second order arithmetic.
We also show that for many natural computable metric spaces and
computable domains the first order theory of the lattice of
effectively open sets is  undecidable. Moreover, for several
important spaces (e.g., $\mathbb{R}^n$, $n\geq1$, and the domain
$P\omega$) this theory is $m$-equivalent to first order
arithmetic.
 \end{abstract}

 \maketitle

\section{Introduction}\label{in}

From  the very beginning of Model Theory, the study of
(un)decidability of first order theories became a central and
popular topic. As a result,  for virtually all structures
$\mathbb{A}=(A;\sigma)$ of a given signature $\sigma$ (in the
Russian literature structures are also known as algebraic systems)
occurring naturally in mathematics their first
order theories $FO(\mathbb{A})$ has been shown to be decidable or undecidable (among vast literature
on the subject we mention \cite{tmr53,eltt,er80}, as examples).

More  recently, several researchers in Computability Theory have been
working on the problem of characterizing the algorithmic complexity
of undecidable first order theories (the complexity is usually
measured by the $m$-degree \cite{ro67,so87} of the theory that in
this context coincides with the 1-1-degree, i.e. with the type of computable isomorphism of
the theory). For many natural structures $\mathbb{A}$ with
undecidable theory  the theory $FO(\mathbb{A})$ turns out to be $m$-equivalent either to first order arithmetic $FO(\mathbb{N})$  or
to  second order arithmetic $FO(\mathbb{N}_2)$ where
$\mathbb{N}:=(\omega;+,\times)$ and $\mathbb{N}_2:=(\omega\cup
P(\omega);\omega,P(\omega),\in,+,\times)$, see e.g.
\cite{ns80,nss96,ni98}.  As is well known (see e.g. \cite{ro67}),
$FO(\mathbb{N})$  is $m$-equivalent to the $\omega$'th iteration
$\emptyset^{(\omega)}$ of the Turing jump  starting from the empty set.

Decidability issues for topological spaces seem to have been studied
less  systematically than for structures arising in algebra, logic
and discrete mathematics, probably because first order language is
not well suited for topology. Nevertheless, there was some important
work for structures related to a topological space $X$, the most
natural of which is the lattice $\mathbf{\Sigma}^0_1(X)$ of open
sets in $X$. To our knowledge, A. Grzegorczyk \cite{gr51} was the
first to consider decidability issues in topology. One of the results in
\cite{gr51} (Corollary 2) interprets  first order arithmetic in
$\mathbf{\Sigma}^0_1(\mathbb{R}^n)$ for each $n\geq2$ which implies
that $FO(\mathbb{N})$ is $m$-reducible to
$FO(\mathbf{\Sigma}^0_1(\mathbb{R}^n))$ and hence the latter theory is
undecidable. The  question of whether
$FO(\mathbf{\Sigma}^0_1(\mathbb{R}))$ is decidable was left open. M.
Rabin \cite{ra69}  answered the question  affirmatively (as well
as the analogous question for the Cantor and Baire spaces) as a corollary
of his  result on decidability of the monadic second order
theory of the binary tree.

A systematic model-theoretic study of structures arising in a topological setting was undertaken in \cite{hjrt77} where it is shown, in particular, that $FO(\mathbb{N}_2)\leq_mFO(\mathbf{\Sigma}^0_1(X))$ for many Hausdorff spaces $X$, and the above mentioned Grzegorczyk's estimate was improved to $FO(\mathbb{N}_2)\equiv_mFO(\mathbf{\Sigma}^0_1(\mathbb{R}^n))$ for each $n\geq2$. Note that in fact the papers \cite{gr51,hjrt77} work with the lattices of closed sets rather than with the lattices of opens but for our purposes this is clearly equivalent.

Another  facet of the relationship between Topology and
Computability Theory is the study of effectivity in a topological
setting as developed in Computable Analysis \cite{wei00} and
Effective Descriptive Set Theory \cite{mo09,s06}. An important
object of study here is the lattices $\Sigma^0_1(X)$ of the so
called effectively open sets in topological spaces $X$ satisfying
some effectivity conditions (see the next section for more details).
The lattice $\Sigma^0_1(X)$ is certainly the most important
sublattice of the lattice $\mathbf{\Sigma}^0_1(X)$ of open sets,
hence it is natural and instructive to study also  definability
and (un)decidability issues for the lattices of effectively open
sets.

This study is interesting and non-trivial even for the discrete
space $\omega$ of natural numbers since in this case the lattice
$\Sigma^0_1(\omega)$ coincides with the lattice $\mathcal{E}$ of
computably  enumerable (c.e.) subsets of $\omega$ which is an
important and popular object of study in Computability Theory
\cite{ro67,so87}. A principal fact about this  lattice is the
undecidability of  $FO(\mathcal{E})$ \cite{he83,he84}. Moreover,
$FO(\mathcal{E})$ is known \cite{hn98} to be $m$-equivalent to 
first-order arithmetic $FO(\mathbb{N})$. Note that the first order
theory of $\mathbf{\Sigma}^0_1(\omega)=P(\omega)$ is decidable
(because $P(\omega)$ is a Boolean algebra).

It seems  that not much is known about (un)decidability of first
order theories of the lattices of effectively open sets except for what is known about the
lattice $\Sigma^0_1(\omega)$ and its relativizations. To our
knowledge, only the cases of the Cantor space $\mathcal{C}$ and the Baire
space $\mathcal{N}$ have been studied to some extent, in the context of
the theory of $\Pi^0_1$-classes. In \cite{ni00} (see the discussion of
Main Theorem 
in the Introduction and Section 3) it is shown that
$FO(\Pi^0_1(\mathcal{C}))$ is $m$-equivalent to $FO(\mathbb{N})$.
Since the lattice $\Pi^0_1(\mathcal{C})$ formed by the complements
of effectively open sets is anti-isomorphic to
$\Sigma^0_1(\mathcal{C})$, this settles the question for
$\mathcal{C}$. To our knowledge,   similar questions for
$\mathcal{N}$ (even the decidability of
$FO(\Sigma^0_1(\mathcal{N}))$) are open.

In this  paper, we make further steps in the study of
(un)decidability issues for the theories of $\Sigma^0_1(X)$ and
$\mathbf{\Sigma}^0_1(X)$. After recalling some necessary preliminaries we
reprove in Section \ref{open} the estimate from \cite{hjrt77}
$FO(\mathbb{N}_2)\equiv_mFO(\mathbf{\Sigma}^0_1(\mathbb{R}^n))$, $n\geq2$,  using the original approach of A. Grzegorczyk which is different from the  approach in \cite{hjrt77}. We also establish the same estimate  for some natural domains. In Section  \ref{efopen} we first
show that for many natural effective spaces $X$ (including 
computable metric spaces without isolated points and many natural
computable domains) the theory $FO(\Sigma^0_1(X))$ is undecidable.
Then we show that $FO(\Sigma^0_1(\mathbb{R}^n))$, $n\geq1$, is
$m$-equivalent to  first order arithmetic. The same
estimate also holds for some natural domains. We conclude in Section
\ref{con} with a discussion of remaining open questions.

The methods  of this paper apply mainly to second countable locally
compact spaces. A precise estimate of the complexity of
$FO(\mathbf{\Sigma}^0_1(X))$ and $FO(\Sigma^0_1(X))$ turns out to be
 subtle and  depends strongly on the topology of $X$. For many
natural spaces $X$ we still have a big gap between the known 
lower and  upper bounds for $FO(\mathbf{\Sigma}^0_1(X))$ and
$FO(\Sigma^0_1(X))$. In particular, for the  Baire space  we currently
only know  the estimate
$\emptyset^\prime\leq_mFO(\Sigma^0_1(X))\leq_mO^{(\omega)}$ where
$O$ is the Kleene ordinal notation system which is a
$\Pi^1_1$-complete set.

Our  upper bounds for the $m$-degree of $FO(\mathbf{\Sigma}^0_1(X))$
and $FO(\Sigma^0_1(X))$ are obtained by a straightforward
application of the Tarski-Kuratowski algorithm, while the lower
bounds use a suitable interpretation of one of the structures
$\mathcal{E},\mathbb{N},\mathbb{N}_2$ in the lattice under
consideration.

This paper is an extended version of the conference paper
\cite{ks16} that contains, in particular,  new results on the lattices of all open sets and an essentially modified proof of Theorem \ref{euclide}  for  $n>1$.

\section{Preliminaries}\label{spaces}

Here we briefly recall some notions and notation relevant to this paper.
We freely  use the standard set-theoretic notation like $|X|$ for
the cardinality of $X$, $X\times Y$ for the cartesian product of
sets and  topological spaces, $P(X)$ for
the set of all subsets of $X$, $\overline{A}$ for the complement $X\setminus A$ of a subset $A$ of a space $X$. 

We assume the reader to be familiar with basic notions of topology
(see e.g. \cite{en89}). We often abbreviate ``topological space'' by
``space''. By $Cl(S)$ (resp. $Int(S)$) we denote the closure (resp. the interior) of a set $S\subseteq X$ in a space $X$. A space $X$ is {\em Polish} if it is separable and
metrizable with a metric $d$ such that $(X,d)$ is a complete metric
space. We denote the set of open subsets of a space $X$ by
$\mathbf{\Sigma}^0_1(X)$. This is the first among the finite levels $\{\mathbf{\Sigma}^0_n(X)\}$ of the Borel hierarchy \cite{ke95,br13} which is formed by applying the
operations of complementation and countable union to the open sets.

Let $\omega$ be the space of non-negative integers with the
discrete topology. The space
$\omega\times\omega=\omega^2$ is
homeomorphic to $\omega$, the  homeomorphism being realized by
the Cantor pairing function $\langle x,y\rangle$.

Let $\omega^\omega$ be the set of all infinite
sequences of natural numbers (i.e., of all functions
$\xi:\omega\rightarrow\omega$). Let $\omega^{<\omega}$ be the set of finite
sequences of elements of $\omega$, including the empty sequence. For
$\sigma\in\omega^{<\omega}$ and $\xi\in{\mathcal N}$, we write
$\sigma\sqsubseteq \xi$ to denote that $\sigma$ is an initial
segment of the sequence $\xi$. By $\sigma\xi=\sigma\cdot\xi$ we
denote the concatenation of $\sigma$ and $\xi$, and by
$\sigma\cdot\mathcal{N}$  the set of all extensions of $\sigma$ in
$\mathcal{N}$. For $x\in\omega^\omega$, we can write
$x=x(0)x(1)\cdots$ where $x(i)\in\omega$ for each $i<\omega$. For
$x\in\mathcal{N}$ and $n<\omega$, let $x[n]=x(0)\cdots x(n-1)$
denote the initial segment of $x$ of length $n$.  Notations in the
style of regular expressions like $0^\omega$, $0^{<\omega} 1$ or $0^m1^n$
have the obvious standard meaning. Define the topology on
$\omega^\omega$ by taking arbitrary unions of sets of the form
$\sigma\cdot\omega^\omega$, where $\sigma\in\omega^{<\omega}$, as the open
sets. The space $\mathcal{N}=\omega^\omega$ with this topology known as the {\em
Baire space}, is of primary importance for Descriptive Set Theory and
Computable Analysis.

For any finite alphabet $A$ (usually  we assume without loss of
generality that $A=k=\{0,\ldots,k-1\}$ where $0<k<\omega$), let
$A^\omega$ be the set of $\omega$-words over $A$. This set may be
topologized similarly to the Baire space. The resulting spaces, which for $k\geq2$ are all (computably) homeomorphic among themselves, are known
as {\em Cantor spaces} (usually the term Cantor space is applied to the
space $\mathcal{C}=2^\omega$ of infinite binary sequences). Note
that the Cantor space is compact while the Baire space is not.

Next we recall  some  definitions related to domain theory (for more
details see e.g. \cite{aj,er93,gh03}). 

Let $X$ be a $T_0$-space. For
$x,y\in X$, let $x\leq y$ denote that $x\in U$ implies $y\in U$, for
all open sets $U$. The relation $\leq$ is a partial order known as
the {\em specialization order}. Let $F(X)$ be the set of {\em
finitary elements} of $X$ (known also as {\em compact elements}),
i.e. elements $p\in X$ such that the upper cone
$\orduparrow{p}=\{x\mid p\leq x\}$ is open. Such open cones are called {\em
$f$-sets}. The space $X$ is called a {\em $\varphi$-space} if every
open set is a union of $f$-sets. Note that every non-discrete
$\varphi$-space is not Hausdorff. A $\varphi$-space is {\em complete} if any non-empty directed set has a supremum w.r.t. the specialization order.

A $\varphi$-space $X$ is an {\em
$f$-space} if any compatible elements $c,d\in F(X)$ have a least
upper bound w.r.t. $\leq$ (compatibility means that $c,d$ have an
upper bound in $F(X)$). An $f$-space $X$ is an {\em $f_0$-space} if
$F(X)$ has a least element.

Let  $\omega^{\leq\omega}=\omega^{<\omega}\cup\omega^\omega$ be the set of
all finite and infinite strings of natural numbers with the topology generated by the sets $\{x\mid u\sqsubseteq x\}$ where $u\in\omega^{<\omega}$. For every $1\leq
k<\omega$, the space $k^{\leq\omega}$ of all finite and infinite words over the alphabet
$\{0,\ldots,k-1\}$ is defined in the same way. Let $P\omega$ be the  powerset of $\omega$ with
the  topology generated by the sets $\orduparrow{F}:=\{A\mid F\subseteq A\subseteq\omega\}$ where $F$ is a finite subset of $\omega$. 

Let $\omega_\bot=\omega\cup\{\bot\}$ be the the space with the topology generated by $\{n\}$, $n\in\omega$. Let
$\omega_\bot^\omega$ be the space of partial functions on
$\omega$ with the topology generated by the sets $\{g\mid f\subseteq g\}$ where $f$ is a function with finite graph and $\subseteq$ is the subgraph relation (as  usual, we identify a partial function $g$ on $\omega$ with the
total function
$\tilde{g}:\omega\rightarrow\omega_\bot=\omega\cup\{\bot\}$ where
$g(x)$ is undefined iff  $\tilde{g}(x)=\bot$). For each $k$, $2\leq
k<\omega$, let $k^\omega_\bot$ be the space of partial functions
$g:\omega\rightharpoonup\{0,\ldots,k-1\}$ defined similarly to
$\omega_\bot^\omega$. 

As is well known,
$\omega^{\leq\omega},k^{\leq\omega},P\omega,\omega_\bot,\omega_\bot^\omega,k^\omega_\bot$ 
are complete $f_0$-spaces where the sets of $f$-elements are respectively $\omega^{<\omega},k^{<\omega}$,the finite subsets of $\omega$, $\omega_\bot$, the finite partial functions on $\omega$, the finite partial functions from $\omega$ to $\{0,\ldots,k-1\}$.

As is well known (see e.g. \cite{er72}), for any
(complete) $f_0$-spaces $X,Y$ the space $Y^X$ of continuous
functions from $X$ to $Y$ with the topology of pointwise convergence
is again a (complete) $f_0$-space. Therefore, any space of continuous partial
functionals over $\omega$ of a finite type is a complete $f_0$-space. In
particular, this applies to the spaces $\mathbb{F}_n$ defined by
induction as $\mathbb{F}_0:=\omega_\bot$,
$\mathbb{F}_{n+1}:=\omega_\bot^{\mathbb{F}_n}$.

Next we explain what we mean by effectively open sets. For any
countably based topological space $X$ and any numbering
$\beta$  of a base of  $X$, define a function
$\pi:\omega\to P(X)$ by $\pi(n)=\bigcup\beta[W_n]$ where $\{W_n\}$
is the standard numbering of the c.e. sets \cite{ro67,so87}) and
$\beta[W_n]=\{\beta(a)\mid a\in W_n\}$. The sets in
$\pi[\omega]$ are called {\em effectively open sets in $X$}. Thus,
the set of  effectively open sets $\Sigma^0_1(X)$ is always equipped
with the induced numbering $\pi$, hence it makes sense to speak
about computable sequences of effectively open sets.

For  many reasonable spaces $X$ with effectivity conditions, one can
define in a natural way the (finite levels of the) effective Borel hierarchy
$\{\Sigma^0_n(X)\}$ and the effective Luzin hierarchy
$\{\Sigma^1_n(X)\}$ (see e.g. \cite{mo09,s06,s15} for details) which
are reasonable effective versions of the classical Borel and Luzin
hierarchies. (Note that the definition of effective hierarchies in a given effective space depends on the chosen numbering of a base of the space.) In particular, $\{\Sigma^0_n(\omega)\}$ and
$\{\Sigma^1_n(\omega)\}$ (taken with a natural numbering of a base in $\omega$) coincide with the  arithmetical and analytical hierarchies of subsets of $\omega$ which are central
objects of study in Computability Theory \cite{ro67}.

We  define some particular classes of effective spaces relevant to
this paper. A {\em computable metric space}  \cite{wei00} is a
triple $(X,d,\nu)$, where $(X,d)$ is a metric space and
$\nu:\omega\rightarrow X$ is a numbering of a dense subset
$rng(\nu)$  of $X$ such that the set
 $$\{(i,j,k,l)\mid\varkappa_k<d(\nu(i),\nu(j))<\varkappa_l\}$$
  is
c.e. Here $\varkappa$ is the conventional numbering of the set
$\mathbb{Q}$ of rationals. Any computable metric space $(X,d,\nu)$
gives rise to a numbering $\beta$ of the standard base
$\beta_{\langle m,n \rangle}=B(\nu_m,\varkappa_n)$ where $\langle m,n
\rangle$ is the Cantor pairing and $B(\nu_m,\varkappa_n)$ is the
basic open ball with center $\nu_m$ and radius $\varkappa_n$ (if $\varkappa_n\leq0$ the ``ball'' is empty).

By a {\em strongly computable metric space}  (SCMS) we mean a
computable metric space such that there exists an infinite
computable sequence $\{B_n\}$ of pairwise disjoint basic open balls.
The metric spaces $\omega$, $\mathbb{Q}$, $\mathcal{C}$,
$\mathcal{N}$, $\mathbb{R}^n$ (where
$\mathbb{R}$ is the space of real numbers) equipped with the
standard metrics and with natural numberings of dense subsets are
SCMS. Any computable metric space without isolated points is an
SCMS. Working with the Euclidean spaces $\mathbb{R}^n$, we denote by $d$ the Euclidean metric, by 0 the zero-vector $(0,\ldots,0)\in\mathbb{R}^n$,  non-empty open (resp. closed) rational balls by $B(a,r)$ (resp. $C(a,r)$) where $a\in\mathbb{Q}^n,r\in\mathbb{Q}^+$. Sometimes it is convenient to use also the empty ball $B(a,0)$.

By a  {\em computable $\varphi$-space} we mean a pair $(X,\delta)$
consisting of a  $\varphi$-space $X$ and a numbering
$\delta:\omega\rightarrow F(X)$ of all the finitary elements such
that the specialization order  is c.e. on the finitary elements
(i.e., the relation $\delta_x\leq\delta_y$ is c.e.). Setting
$\beta(n):=\orduparrow{\delta_n}$ we obtain a numbering of a
topological base of $X$. Thus, we have a notion of an effective
open set in every computable $\varphi$-space.

By a  {\em strongly computable $\varphi$-space} (SC$\Phi$S) we mean a
computable $\varphi$-space  $X$ such that the specialization order
is computable on the finitary elements, and there is a computable
sequence $\{c_n\}$ of pairwise incomparable finitary elements.
 An SC$\Phi$S $X$ is a
{\em strongly computable $f_0$-space (SCF$_0$S)} if it is an
$f_0$-space, the relation of compatibility is computable on $F(X)$,
and the supremum of compatible finitary elements is computable.
Although the restrictions imposed on SC$\Phi$Ss and SCF$_0$Ss are rather strong, many popular
domains are SCF$_0$Ss. In particular, this applies to all concrete
examples of $\varphi$-spaces mentioned above in this section. (For
the space $Y^X$ of continuous functions, a close inspection of the
corresponding proofs \cite{er72} shows that if $X,Y$ are (complete)
SCF$_0$Ss then so is  $Y^X$. Therefore, any space of continuous
partial functionals over $\omega$ of a finite type is a complete
SCF$_0$S.) Note that the ``strong'' variations above are rather ad hoc and do not pretend to be fundamental notions in the field.

We  conclude this section by briefly recalling of some notions from logic. We
consider only structures of finite relational signatures (when a
functional symbol is used, as e.g. in the  structure $\mathbb{N}$,
we identify the corresponding function with its graph). For a
$\sigma$-structure $\mathbb{A}=(A;\sigma)$, a relation $R\subseteq
A^k$ is {\em definable in $\mathbb{A}$} if there is
a first order $\sigma$-formula $\phi(x_1,\ldots,x_k,p_1,\ldots,p_l)$
and (possibly) some values $p_1,\ldots,p_l\in A$ of parameters such that
 $$R=\{(x_1,\ldots,x_k)\in A^k\mid\mathbb{A}\models \phi(x_1,\ldots,x_k,p_1,\ldots,p_l)\}.$$
 If the  list of
parameters is empty then we speak about {\em definability without parameters}. Thus, $R\subseteq A^k$ is {\em definable in $\mathbb{A}$
without parameters} if there is a first order $\sigma$-formula
$\phi(x_1,\ldots,x_k)$ with
 $$R=\{(x_1,\ldots,x_k)\in
A^k\mid\mathbb{A}\models \phi(x_1,\ldots,x_k)\}.$$
  A function on $A$ is definable (with or without parameters)  if its graph is definable. An
element of $A$ is definable if the corresponding singleton set
$\{a\}$ is definable. A structure is definable if its universe and
all signature predicates are definable.

E.g., if $\mathbb{A}=(A;\cup,\cap,0,1)$ is a bounded distributive
lattice then any of $\cup,\cap,0,1$ is definable without parameters
in $(A;\leq)$ where $\leq$ is the induced partial order on $A$, and
$\leq$ is definable in $(A,\cup)$. Moreover, in this case we can
even speak about arbitrary Boolean terms of elements of $A$ (meaning
their values in a Boolean algebra extending $\mathbb{A}$). Thus,
dealing with our lattices $\mathbf{\Sigma}^0_1(X)$ and
$\Sigma^0_1(X)$ we can mean any of the signatures $\{\subseteq\}$,
$\{\cup,\cap,\emptyset,X\}$, or even $\{\cup,\cap,\bar{
},\emptyset,X\}$. For simplicity, we omit the signature symbols in
the notation of these structures. We use in our formulas some standard abbreviations, in particular the bounded quantifiers $(\forall x)_{\phi(x)}\psi:=\forall x(\phi(x)\to\psi)$ and $(\exists x)_{\phi(x)}\psi:=\exists x(\phi(x)\wedge\psi)$ or the ``quantifier'' $\exists! x\psi$ meaning ``there exists a unique $x$ satisfying $\psi$''.

The  first order theory $FO(\mathbb{A})$ of the structure
$\mathbb{A}$ is the set of $\sigma$-sentences true in $\mathbb{A}$.
Along with first order logic, in logical theories some other logics
are considered, in particular the (monadic) second order logic where
one can use, along with the usual variables, also variables ranging
over the (unary) relations on $A$. Accordingly, one can consider the
monadic second order theory $\mathit{MSO}(\mathbb{A})$, or the full second order
theory $SO(\mathbb{A})$, of $\mathbb{A}$ (in the latter case one needs
variables for relations of any arity). Since the Cantor coding is
definable in $\mathbb{N}$ without parameters,
$SO(\mathbb{N})\equiv_m\mathit{MSO}(\mathbb{N})$. The theory
$\mathit{MSO}(\mathbb{A})$ may be considered as the first order theory
$FO(\mathbb{A}^{mso})$ of the extended structure
 $$\mathbb{A}^{mso}:=(A\cup P(A);A,P(A),\in,\sigma)$$
  obtained from
$\mathbb{A}$ by adjoining the powerset of $A$ to the universe, the
unary predicates for $A$ and $P(A)$, and the membership relation to
the list of relations. In particular, we have
$\mathbb{N}_2=\mathbb{N}^{mso}$ and
$\mathit{SO}(\mathbb{N})\equiv_mFO(\mathbb{N}_2)$.

An  important tool to compare algorithmic complexity of theories is
the notion of interpretability of one theory or structure in
another. In fact, there are many versions of this notion (see e.g.
\cite{tmr53,eltt,er80} of which we briefly recall  a couple of those
used in the sequel.

A  $\tau$-structure $\mathbb{B}$ is {\em interpretable in a
$\sigma$-structure $\mathbb{A}$ without parameters} if some
isomorphic copy of $\mathbb{B}$ is definable in $\mathbb{A}$ without
parameters. A weaker version of this is the notion of
$c$-interpretability (where $c$ comes from ``congruence''). We say
that $\mathbb{B}$ is  $c$-interpretable in $\mathbb{A}$ without
parameters if there exist a $\tau$-structure $\mathbb{C}$ and a
congruence $\sim$ on $\mathbb{C}$ such that both $\mathbb{C}$ and
$\sim$ are definable in $\mathbb{A}$ without parameters and the quotient-structure
$\mathbb{C}/\mathord\sim$ (whose elements are the equivalence classes $c/\mathord\sim,c\in C$) is isomorphic to $\mathbb{B}$. Interpretability
with parameters is introduced in the same manner.

As is well known (see e.g. \cite{tmr53,eltt}), if $\mathbb{B}$   is $c$-interpretable in $\mathbb{A}$
without parameters   then  $FO(\mathbb{B})\leq_mFO(\mathbb{A})$. The
same is true for definability with parameters provided that the set
of ``defining'' parameters may be chosen definable. The latter notion means that
there is a non-empty definable (without parameters) set $P\subseteq
A^l$ of parameters such that for any value of parameters in $P$ the
corresponding structure is isomorphic to $\mathbb{B}$. For  future reference, we formulate some of the mentioned facts as a lemma. 

\begin{lem}\label{inter}
Let $\mathbb{B}$   be $c$-interpretable in $\mathbb{A}$
without parameters (or with a non-empty set of parameters which is itself definable without parameters).  Then  $FO(\mathbb{B})\leq_mFO(\mathbb{A})$.
 \end{lem}

When no non-empty set of eligible parameters is definable, the relation
$FO(\mathbb{B})\leq_mFO(\mathbb{A})$ is not true in general  but
there is some version of undecidability which is preserved also by
such interpretations. A theory (not necessarily complete) of
signature $\sigma$ is  {\em hereditarily undecidable} if any of its
subtheories of signature $\sigma$ is undecidable. It is well known
(see e.g. \cite{er80}) that if $FO(\mathbb{B})$ is hereditarily
undecidable and $\mathbb{B}$   is $c$-interpretable in $\mathbb{A}$
with parameters then $FO(\mathbb{A})$ is hereditarily undecidable.

Additional information about interpretations may be found on page 215 of \cite{ho93}.

\section{The lattices of open sets}\label{open}

Here we give precise estimates of the algorithmic complexity of
$FO(\mathbf{\Sigma}^0_1(X))$ for some spaces~$X$.

\subsection{An upper bound}

First we establish a natural upper bound that applies to many
countably based locally compact spaces. We need a technical notion
related to local compactness. By {\em analytically locally compact
space} (AnLCS) we mean a triple $(X,\beta,\kappa)$ consisting of a
topological space $X$, a numbering $\beta$ of a base in $X$
containing  the empty set (the presence of the empty set is not principal and maybe removed by using slight modification of the notion of AnLCS), and a numbering $\kappa$ of some compact
sets in $X$ such that any set $\beta_n$ is a union of some sets in
$\{\kappa_i\mid i<\omega\}$, and the relation
$\kappa_i\subseteq\bigcup\beta[D_n]$ (where $\{D_n\}$
is the canonical numbering of finite subsets of $\omega$
\cite{ro67}) is analytical, i.e. it is in $\bigcup_n\Sigma^1_n(\omega)$.

Note that although AnLCSs are not automatically locally compact,
many locally compact spaces may be considered as AnLCSs. In
particular, the computable $\varphi$-spaces, the finite dimensional
Euclidean spaces, and the Cantor space are AnLCSs (for instance, for a computable
$\varphi$-space $(X,\delta)$ we can set
$\kappa_n:=\beta_n:=\orduparrow{\delta_n}$ which is compact;  the relation
$\kappa_i\subseteq\bigcup\beta[D_n]$ in this case is c.e.).

\begin{prop}\label{loc2}
If $(X,\beta,\kappa)$ is an AnLCS then $FO(\mathbf{\Sigma}^0_1(X))\leq_mFO(\mathbb{N}_2)$.
\end{prop}

\proof  Define a surjection
$\tau:\mathcal{N}\to\mathbf{\Sigma}^0_1(X)$  by
$\tau(p)=\bigcup_n\beta_{p(n)}$. The surjection $\tau$ has several
nice properties, in particular it is an admissible representation of
the hyperspace of open sets in $X$ (see e.g. \cite{s13} for
additional details).

It suffices  to show that the relation $\tau(p)\subseteq\tau(q)$ is
an analytical subset of $\mathcal{N}\times\mathcal{N}$ because then
the elementary diagram of the represented structure
$(\mathbf{\Sigma}^0_1(X);\subseteq,\tau)$, and hence also
$FO(\mathbf{\Sigma}^0_1(X))$, are $m$-reducible to
$FO(\mathbb{N}_2)$.

Obviously,  $\tau(p)\subseteq\tau(q)$ is equivalent to $\forall
n(\kappa_n\subseteq\tau(p)\rightarrow\kappa_n\subseteq\tau(q))$,
hence it suffices to show that the relation
$\kappa_n\subseteq\tau(p)$ is analytical. We  have
$\kappa_n\subseteq\tau(p)$ iff
$\kappa_n\subseteq\bigcup_i\beta_{p(i)}$ iff $\exists
m(\kappa_n\subseteq\beta_{p(0)}\cup\cdots\cup\beta_{p(m)})$, by
compactness of $\kappa_n$. The set
 $$\{(m,n,p)\mid\kappa_n\subseteq\beta_{p(0)}\cup\cdots\cup\beta_{p(m)}\})$$
is, by the definition of AnLCS, an analytical subset of
$\omega\times\omega\times\mathcal{N}$, hence the relation
$\kappa_n\subseteq\tau(p)$ is analytical.
 \qed

\subsection{Lattices of opens in Euclidean spaces}

The next result improves the estimate from  \cite{gr51} mentioned in
the Introduction, providing a different proof of this result compared with \cite{hjrt77}.

\begin{thm}\label{euclide2}
For any $n\geq2$,
$FO(\mathbf{\Sigma}^0_1(\mathbb{R}^n))\equiv_mFO(\mathbb{N}_2)$.
\end{thm}

\proof Since the upper bound holds by 
Proposition \ref{loc2}, we only have to prove the lower  bound
$FO(\mathbb{N}_2)\leq_mFO(\mathbf{\Sigma}^0_1(\mathbb{R}^n))$. For
this, we extend  Grzegorczyk's interpretation \cite{gr51} of
$\mathbb{N}$ in $\mathbf{\Sigma}^0_1(\mathbb{R}^n)$ to an
interpretation  of $\mathbb{N}_2$ in
$\mathbf{\Sigma}^0_1(\mathbb{R}^n)$. We start with a brief sketch of
the Grzegorczyk interpretation (with slightly different notation).

Since $\mathbf{\Sigma}^0_1(\mathbb{R}^n)$ is a distributive lattice,
we can use in the definitions not only the symbol of inclusion but
also the symbols of Boolean operations and the constants
$\emptyset,\mathbb{R}^n$. Note that for any $x\in \mathbb{R}^n$ the
set $\mathbb{R}^n\setminus\{x\}$ is open, and the class of such
co-singleton sets is definable  in
$\mathbf{\Sigma}^0_1(\mathbb{R}^n)$ as the class of sets maximal
w.r.t. inclusion among the sets strictly below $\mathbb{R}^n$. We also use some other observations from \cite{hjrt77}.

The
key observation of Grzegorczyk was the definability in
$\mathbf{\Sigma}^0_1(\mathbb{R}^n)$ (without parameters) of the set
$\mathit{Cof}$ of cofinite subsets of $\mathbb{R}^n$, as
well as of the following relations on $\mathit{Cof}$:
\begin{itemize}
\item $U\approx V$ iff $|\overline{U}|=|\overline{V}|$ where
$\overline{U}:=\mathbb{R}^n\setminus U$,

\item $P_+(U,V,W)$ iff $|\overline{U}|+|\overline{V}|=|\overline{W}|$,

\item $P_\times(U,V,W)$ iff
$|\overline{U}|\times|\overline{V}|=|\overline{W}|$.
\end{itemize}

Grzegorczyk's interpretation is now given by the natural
isomorphism $k\mapsto U_k/\mathord\approx$ from $\mathbb{N}$ onto the quotient
of $(\mathit{Cof};P_+,P_\times)$ modulo the definable congruence $\approx$ where
$U_k:=\mathbb{R}^n\setminus\{(i,\ldots,i)\mid i<k\}$.

The definability of $\mathit{Cof},\approx,P_+,P_\times$ makes  heavy use of the
notions of connected open sets and open connected components which are
definable, respectively, by the formulas
	 $$\mathit{Con}(U):=U\not=\emptyset\wedge\neg\exists V,W(U\subseteq V\cup
W\wedge V\cap W=\emptyset\wedge U\cap
  V\not=\emptyset\wedge U\cap W\not=\emptyset)$$
 and
	 $$\mathit{Cmp}(U,V):=U\subseteq V\wedge \mathit{Con}(U)\wedge\forall U^\prime(U^\prime\subseteq V\wedge
	 \mathit{Con}(U^\prime)\wedge U\cap U^\prime\not=\emptyset\rightarrow U^\prime\subseteq U).$$
 Let now $\xi(U,V)$ be a formula saying in
$\mathbf{\Sigma}^0_1(\mathbb{R}^n)$ that $U$ is coinfinite, $V$ is
cofinite, $U\cap F=\emptyset$ (where $F$ is the finite set of all  points in the complement of $V$), $U\cup F$
is open, and $V$ is the smallest element of $(\mathit{Cof};\subseteq)$ with
these properties. In other words,
$\mathbf{\Sigma}^0_1(\mathbb{R}^n)\models\xi(U,V)$ iff $U$ is
obtained from the open set $U\cup F$ by removing the
points from $F$; Note that, by the minimality of $V$, for any
coinfinite open set $U$ there is at most one cofinite open set $V$
with $\mathbf{\Sigma}^0_1(\mathbb{R}^n)\models\xi(U,V)$.

Let $\mathcal{P}$ be the definable subset of
$\mathbf{\Sigma}^0_1(\mathbb{R}^n)$ formed by the coinfinite sets
$W$ satisfying $(\forall U)_{\mathit{Cmp}(U,W)}\exists V\xi(U,V)$. Then
$\mathcal{P}$ is disjoint with $\mathit{Cof}$ and we can associate with any
$W\in\mathcal{P}$ the set
 $$A_W:=\{k\in\omega\mid(\exists U)_{\mathit{Cmp}(U,W)}\exists V(\xi(U,V)\wedge|\overline{V}|=k)\}.$$

Note that for any $A\subseteq\omega$ there is $W=W(A)\in\mathcal{P}$
with $A=A_W$. Indeed, we can take the set $W:=\bigcup\{S_a\mid a\in
A\}$ where $S_a$ is obtained from the open set
$(a,a+1)\times\mathbb{R}^{n-1}$ by removing $a$ points (note that
$S_a$ are the connected components of $W$). Note also that the
relations
 $$\epsilon(V,W)\Leftrightarrow V\in \mathit{Cof}\wedge
W\in\mathcal{P}\wedge|\overline{V}|\in A_W$$
 and
 $$W\equiv
W^\prime\Leftrightarrow W\in\mathcal{P}\wedge
W^\prime\in\mathcal{P}\wedge A_W=A_{W^\prime}$$
 are definable in
$\mathbf{\Sigma}^0_1(\mathbb{R}^n)$ without parameters.

These remarks mean that the maps $k\mapsto U_k/\mathord\approx$, $A\mapsto
W(A)/\mathord\equiv$ give an isomorphism  from $\mathbb{N}_2$ onto the
quotient of
$(\mathit{Cof}\cup\mathcal{P};\mathit{Cof},\mathcal{P},\epsilon,P_+,P_\times)$ modulo
the definable congruence on $\mathit{Cof}\cup\mathcal{P}$ induced by
$\approx$ and $\equiv$. Thus, $\mathbb{N}_2$ is $c$-interpretable in
$\mathbf{\Sigma}^0_1(\mathbb{R}^n)$ without parameters. By Lemma \ref{inter},
$FO(\mathbb{N}_2)\leq_mFO(\mathbf{\Sigma}^0_1(\mathbb{R}^n))$.
 \qed

\subsection{Lattices of opens in domains}

Here we give  a similar estimate for some natural domains.

\begin{thm}\label{scott2}
 For any $X\in\{P\omega,k^\omega_\bot\mid2\leq k\leq\omega\}$,
 $FO(\mathbf{\Sigma}^0_1(X))\equiv_mFO(\mathbb{N}_2)$.
\end{thm}

\proof  First consider the  space $P\omega$ (which is homeomorphic to $1^\omega_\bot$). Since the upper bound
holds by Proposition \ref{loc2}, we only have to prove the lower
bound $FO(\mathbb{N}_2)\leq_mFO(\mathbf{\Sigma}^0_1(P\omega))$. By an observation of R. Robinson mentioned in the beginning of the proof of Lemma 7.2 in \cite{hjrt77}, $\times$ is definable in $(\omega;+)$ by a monadic  second order formula, so if suffices to prove  $FO(\mathbb{N}^\prime_2)\leq_mFO(\mathbf{\Sigma}^0_1(P\omega))$ where $\mathbb{N}^\prime_2:=(\omega\cup P(\omega);\omega,P(\omega),\in, +)$.

We show that the class $\mathcal{F}$ of $f$-sets (i.e., sets of the form
$\orduparrow{F}$ where $F\subseteq\omega$ is a finite set),  is definable in
$\mathbf{\Sigma}^0_1(P\omega)$ without parameters. Indeed, a
defining formula is
 $$ir(V):=V\not=\emptyset\wedge\forall U,U^\prime(V\subseteq
U\cup U^\prime\rightarrow V\subseteq U\vee V\subseteq U^\prime)$$
 which says in $\mathbf{\Sigma}^0_1(P\omega)$ that
$V$ is a non-zero join-irreducible element. Obviously, any set
$V=\orduparrow{F}$ satisfies this formula. Conversely, let an element
$V\in\mathbf{\Sigma}^0_1(P\omega)$ satisfy the formula. Since
$V\not=\emptyset$, for some  sequence $\{F_n\}$ of finite sets we
have $V=\bigcup_n\orduparrow{F}_n$. Since the partial order $(\{F_n\mid
n<\omega\};\subseteq)$ is well founded, it has a minimal element
$F$. Then of course $\orduparrow{F}\subseteq V$, so it suffices to show
that $V\subseteq\orduparrow{F}$.
Let $S:=\bigcup\{\orduparrow{F_n}\mid F\not\subseteq F_n\}$. Since
$S\in\mathbf{\Sigma}^0_1(P\omega)$, $V\subseteq\orduparrow{F}\cup S$
and $V$ is join-irreducible, it suffices to show that
$V\not\subseteq S$. Suppose the contrary, then $F\in S$, so
$F_n\subseteq F$ for some $n$ with $F\not\subseteq F_n$,
contradicting the minimality of $F$.

Let  $\mathcal{G}$ be the set of finite subsets of $\omega$. Since
$(\mathcal{G};\subseteq)$ is isomorphic to $(\mathcal{F};\supseteq)$
via $G\mapsto\orduparrow G$, we can  interpret $\mathbb{N}^\prime_2$ in
$\mathbf{\Sigma}^0_1(P\omega)$ similarly to the previous proof. Define
the relations $\approx,P_+$ on $\mathcal{F}$ as follows:

\begin{itemize}
\item $\orduparrow{F}\approx\orduparrow{G}$ iff $|F|=|G|$,

\item $P_+(\orduparrow{F},\orduparrow{G},\orduparrow{H})$ iff $|F|+|G|=|H|$.


\end{itemize}
We show that $\approx,P_+,P_\times$ are definable in $\mathbf{\Sigma}^0_1(P\omega)$ without parameters.

First  we show that the set $\mathcal{V}_n:=\{\orduparrow{F}:|F|=n\}$
is definable in  $\mathbf{\Sigma}^0_1(P\omega)$ for each $n<\omega$.
A sequence $\{\phi_n(V)\}$ of defining formulas is given by
induction on $n$ as follows: $\phi_0(V):=\forall U(U\subseteq V)$
(saying that $V$ is the largest element of a lattice), and
 $$\phi_{n+1}(V):=\neg\phi_0(V)\wedge\cdots\wedge
 \neg\phi_n(V)\wedge( \forall U)_{ir(V)}(V\subset U\rightarrow\phi_0(V)
 \vee\cdots\vee\phi_n(V))$$
(saying in $\mathbf{\Sigma}^0_1(P\omega)$ that $V$ is a maximal
join-irreducible element among those not in the set defined by the
formula $\phi_0(V)\vee\cdots\vee\phi_n(V)$). By induction on $n$,
$\phi_n$ defines $\mathcal{V}_n$ for each $n$.

Note that  any $U\in\mathbf{\Sigma}^0_1(P\omega)$ is uniquely
representable as the union of its  maximal
(w.r.t. inclusion) $f$-subsets. Clearly, the relation
$\mathit{Maxf}(C,U)$ meaning that ``$C$ is a maximal $f$-subset of $U$'' is definable
in $\mathbf{\Sigma}^0_1(P\omega)$.

Let  $A\approx_UB$ be the definable (in
$\mathbf{\Sigma}^0_1(P\omega)$) relation
``$A=\orduparrow{F},B=\orduparrow{G}\in\mathcal{F}$, $F\cap G=\emptyset$,
$|F\cap H|=|G\cap H|=1$ for each maximal $f$-subset $\orduparrow{H}$ of $U$,
$H\cap H^\prime=\emptyset$, for all distinct maximal $f$-subsets
$\orduparrow{H},\orduparrow{H}^\prime$ of $U$, and $F\cup
G\subseteq\bigcup\{H\mid \mathit{Maxf}(\orduparrow{H},U)\}$''. For instance, $\orduparrow{\emptyset}\approx_\emptyset\orduparrow{\emptyset}$. Clearly,
$\orduparrow{F}\approx_U\orduparrow{G}$ implies $|F|=|G|$, and for all
$F,G\in\mathcal{F}$ we have: $|F|=|G|$ iff $|F\setminus
G|=|G\setminus F|$ iff $\exists U(\orduparrow{(F\setminus
G)}\approx_U\orduparrow{(G\setminus F)})$. This yields the definability
of $\approx$.

Since  $|F|+|G|=|H|$ iff there are disjoint
$F^\prime,G^\prime\in\mathcal{G}$ such that $|F^\prime|=|F|$,
$|D^\prime|=|G|$ and $|F^\prime\cup G^\prime|=|H|$, the relation
$P_+$ is definable.


Since  $k\mapsto\orduparrow{\{0,\ldots,k-1\}}/\mathord\approx$ is an
isomorphism from $(\omega;+)$ onto the quotient of
$(\mathcal{F};P_+)$ modulo $\approx$, we get a
$c$-interpretation of $(\omega;+)$ in $\mathbf{\Sigma}^0_1(P\omega)$
without parameters. As in the proof of the previous theorem, this
interpretation can be easily  extended to a $c$-interpretation of
$\mathbb{N}_2$ in $\mathbf{\Sigma}^0_1(P\omega)$ without parameters.

Namely, let
$\mathcal{P}:=\mathbf{\Sigma}^0_1(P\omega)\setminus\mathcal{F}$, so
any $W\in\mathcal{P}$ is either empty or has at least two
maximal $f$-subsets. Then $\mathcal{P}$ is definable and we can associate with any
$W\in\mathcal{P}$ the set
 $$A_W:=\{k-1\in\omega\mid\exists K\in\mathcal{G}(\mathit{Maxf}(\orduparrow{K},W)\wedge k=|K|)\}.   $$
 Note that for any $A\subseteq\omega$ there is $W=W(A)\in\mathcal{P}$ with $A=A_W$.
If $A$ is empty, we have to take $W=\emptyset$. Otherwise, fix
$a_0\in A$ and choose a countable partition $Q,R_0,R_1,\ldots$ of
$\omega$ into infinite sets with $Q=\{q(0)<q(1)<\cdots\}$ and
$R_i=\{r_i(0)<r_i(1)<\cdots\}$; now it suffices to set
 $$W(A):=\orduparrow{\{q(0),\ldots,q(a_0)\}}\cup\bigcup\{\orduparrow{\{r_a(0),\ldots,r_a(a)\}}\mid a\in A\}.$$
Furthermore, the relations
 $$\epsilon(\orduparrow{F},W)\Leftrightarrow F\in\mathcal{G}\wedge
W\in\mathcal{P}\wedge|F|\in A_W$$
 and
 $$W\equiv
W^\prime\Leftrightarrow W\in\mathcal{P}\wedge
W^\prime\in\mathcal{P}\wedge A_W=A_{W^\prime}$$
 are definable in
$\mathbf{\Sigma}^0_1(P\omega)$ without parameters. Therefore, the
maps $k\mapsto\orduparrow{\{0,\ldots,k-1\}}/\mathord\approx$, $A\mapsto
W(A)/\mathord\equiv$ give an isomorphism  from $\mathbb{N}^\prime_2$ onto the
quotient of
$(\mathcal{F}\cup\mathcal{P};\mathcal{F},\mathcal{P},\epsilon,P_+)$
modulo the definable congruence on $\mathcal{F}\cup\mathcal{P}$
induced by $\approx$ and $\equiv$. Thus, $\mathbb{N}^\prime_2$ is
$c$-interpretable in $\mathbf{\Sigma}^0_1(P\omega)$ without
parameters. By Lemma \ref{inter},
$FO(\mathbb{N}^\prime_2)\leq_mFO(\mathbf{\Sigma}^0_1(P\omega))$, completing the proof for $P\omega$.

Now consider the space $k^\omega_\bot$. Since again $FO(\mathbf{\Sigma}^0_1(k^\omega_\bot))\leq_mFO(\mathbb{N}_2)$ by Proposition \ref{loc2}, it suffices to show $FO(\mathbb{N}_2)\leq_mFO(\mathbf{\Sigma}^0_1(k^\omega_\bot))$. By the just established estimate for $P\omega$, it suffices to show $FO(\mathbf{\Sigma}^0_1(P\omega))\leq_mFO(\mathbf{\Sigma}^0_1(k^\omega_\bot))$, i.e., to interpret $\mathbf{\Sigma}^0_1(P\omega)$ in $\mathbf{\Sigma}^0_1(k^\omega_\bot)$ with a non-empty definable set of parameters.

For any total function $f:\omega\to k$, let $A_f$ be the set of partial subfunctions of $f$. Then $A_f$ is a closed subset of $k^\omega_\bot$ which, taken with the subspace topology, is homeomorphic to $1^\omega_\bot$, and hence also to $P\omega$. Therefore, $\mathbf{\Sigma}^0_1(P\omega)$ is isomorphic to $\mathbf{\Sigma}^0_1(A_f)$ for each $f:\omega\to k$, so it suffices to show that $\mathcal{A}:=\{A_f\mid f:\omega\to k\}$ is definable in $\mathbf{\Sigma}^0_1(k^\omega_\bot)$ without parameters. This follows from the following assertion that is easy to check: $A\in\mathcal{A}$ iff $A$ is a maximal (w.r.t. inclusion) closed subset of $k^\omega_\bot$ satisfying
\begin{equation*}
 \forall U,V\in\mathbf{\Sigma}^0_1(k^\omega_\bot)(U\cap A\not=\emptyset\wedge
V\cap A\not=\emptyset\rightarrow U\cap V\cap A\not=\emptyset).
\tag*{\qEd}
\end{equation*}

As mentioned in the Introduction, M. Rabin has shown  that
$FO(\mathbf{\Sigma}^0_1(X))$ is decidable for
$X\in\{\mathcal{C},\mathcal{N},\mathbb{R}\}$. The same holds for
$X\in\{[0,1], (0,1]\}$, with a slight modification of the proof in
\cite{ra69}. We conclude this section with a natural variation of
this for domains.

\begin{prop}\label{dom2}
 For any $1\leq k\leq\omega$, $FO(\mathbf{\Sigma}^0_1(k^{\leq\omega}))$ is decidable.
\end{prop}

\proof Relate to any $A\subseteq k^{<\omega}$ the open  set
$U_A:=\bigcup_{a\in A}a\cdot k^{\leq\omega}$. Then $A\mapsto U_A$ is
a surjection from $P(k^{<\omega})$ onto
$\mathbf{\Sigma}^0_1(k^{\leq\omega})$. Let $\preceq$ be the preorder
on $P(k^{<\omega})$ defined by: $A\preceq B$ iff $U_A\subseteq U_B$. Then
the lattice $\mathbf{\Sigma}^0_1(k^{\leq\omega})$ is isomorphic to
the quotient order of $(P(k^{<\omega});\preceq)$ modulo the induced
congruence $\approx$, hence $\mathbf{\Sigma}^0_1(k^{\leq\omega})$ is
$c$-interpretable in $(P(k^{<\omega});\preceq)$ without parameters and
therefore
$FO(\mathbf{\Sigma}^0_1(k^{\leq\omega}))\leq_mFO(P(k^{<\omega});\preceq)$.

Since $A\preceq B$ iff $\forall a\in A(a\cdot
k^{\leq\omega}\subseteq U_B)$ iff $\forall a\in A(a\in U_B)$ iff
$\forall a\in A\exists b\in B(b\sqsubseteq a)$, the structure
$(P(k^{<\omega});\preceq)$ is interpretable in
	$(k^{<\omega};\sqsubseteq)^{\mathit{mso}}$ without parameters and therefore
$FO(P(k^{<\omega});\preceq)\leq_m \mathit{MSO}(k^{<\omega};\sqsubseteq)$. Since the
last theory is decidable by \cite{ra69},
$FO(\mathbf{\Sigma}^0_1(k^{\leq\omega}))$ is decidable.
 \qed

\section{The lattices of effectively open sets}\label{efopen}

In this section we examine algorithmic complexity of the first order
theories of  lattices of effectively open sets.

\subsection{Preliminary results}

We start by
showing that for many natural effective spaces their theories of the effectively open sets are
hereditarily undecidable.

\begin{thm}\label{metric}
Let $X$  be an SCMS or an SC$\Phi$S. Then $FO(\Sigma^0_1(X))$ is hereditarily undecidable.
\end{thm}

\proof Since  the theory $FO(\mathcal{E})$ is hereditarily
undecidable \cite{he83,he84}, it suffices to  interpret
$\mathcal{E}$ in  $\Sigma^0_1(X)$  with parameters; in our case
two parameters $V,W$ will suffice.

Consider  the formulas $\phi(U,V,W):=V\subseteq U\wedge U\subseteq
W$ and
 $$\phi_\subseteq(U,U^\prime,V,W):=\phi(U,V,W)\wedge\phi(U^\prime,V,W)\wedge
U\subseteq U^\prime.$$
 For any values $V,W\in\Sigma^0_1(X)$ for parameters
with $V\subseteq W$,  the formulas define the sublattice
$\mathcal{D}$ of $\Sigma^0_1(X)$ formed by the sets lying between
$V$ and $W$. Therefore, it suffices to find effectively open sets
$V,W$ such that the lattice $\mathcal{D}$ is isomorphic to
$\mathcal{E}$.

Let first $X$  be an SCMS. Let $\{B_n\}$ be the sequence of basic
open balls from the definition of SCMS and let $B_n^\prime$  be
obtained from the ball $B_n$ by removing its center $c_n$. Then
$\{B_n^\prime\}$ is a computable sequence of effectively open sets,
hence $V:=\bigcup_nB_n^\prime$ and $W:=\bigcup_nB_n$ are effectively
open and $W\setminus V=\{c_0,c_1,\ldots\}$. From the definition of
SCMS it is easy to see that the function $D\mapsto\{n\mid c_n\in
D\}$ is a desired isomorphism between $\mathcal{D}$ and
$\mathcal{E}$.

Now let $X$ be an SC$\Phi$S.  Let $\{c_n\}$ be the sequence of finitary
elements from the definition of SC$\Phi$S.  This time we take as 
parameters the values $W:=\bigcup_n\orduparrow{c_n}$ and
$V:=\bigcup_n(\orduparrow{c_n}\setminus \{c_n\})$. From the definition
of SC$\Phi$S it is easy to see that $\{\orduparrow{c_n}\setminus \{c_n\}\}$
is a computable sequence of effectively open sets, hence again $V$
and $W$ are effectively open and $W\setminus V=\{c_0,c_1,\ldots\}$.
Moreover, the function $D\mapsto\{n\mid c_n\in D\}$ is again an
isomorphism between $\mathcal{D}$ and $\mathcal{E}$.
 \qed

Next we aim to obtain precise estimates of the algorithmic complexity of
$FO(\Sigma^0_1(X))$ for some spaces $X$.
First  we establish a natural upper bound that applies to many
countably based locally compact spaces.  {\em Arithmetically locally
compact spaces} (ArLCS) are defined precisely as AnLCS in the
previous section but this time the relation
$\kappa_i\subseteq\bigcup\beta[D_n]$ is required to be arithmetical, i.e. to be in $\bigcup_n\Sigma^0_n(\omega)$.
Again, many locally compact spaces may be considered as
ArLCSs. In particular, the computable $\varphi$-spaces, the finite
dimensional Euclidean spaces, and the Cantor space are ArLCSs.

The following proposition is an easy variation of Proposition \ref{loc2}.

\begin{prop}\label{loc}
If $(X,\beta,\kappa)$ is an ArLCS then
$FO(\Sigma^0_1(X))\leq_mFO(\mathbb{N})$.
\end{prop}

\proof Recall that $\pi_i:=\bigcup\beta[W_i]$ is the natural numbering of effectively open sets. It suffices  to show that the relation
$\pi_i\subseteq\pi_j$ is arithmetical because then the elementary
diagram of the numbered structure $(\Sigma^0_1(X);\subseteq,\pi)$,
and hence also $FO(\Sigma^0_1(X))$, is $m$-reducible to $FO(\mathbb{N})$.

Obviously,  $\pi_i\subseteq\pi_j$ is equivalent to $\forall
n(\kappa_n\subseteq\pi_i\rightarrow\kappa_n\subseteq\pi_j)$, hence
it suffices to show that the relation $\kappa_n\subseteq\pi_i$ is
arithmetical. We have:
$\kappa_n\subseteq\pi_i$ iff
$\kappa_n\subseteq\bigcup\beta[W_i]$ iff $\exists m(D_m\subseteq W_i
\wedge \kappa_n\subseteq\bigcup\beta[D_m])$
(the last equivalence holds by compactness of
$\kappa_n$). The last relation is arithmetical by the definition of
ArLCS.
 \qed

In the next subsection we give   precise estimates of the $m$-degrees of
$FO(\Sigma^0_1(\mathbb{R}^n))$  for which we need the following lemma. Recall that the definition of $\Sigma^0_1(\mathbb{R}^n))$ depends on the natural numbering of rational open balls in $\mathbb{R}^n)$.
 
\begin{lem}\label{con}
 Any connected  component of an effectively open set in $\mathbb{R}^n$ is effectively open.
\end{lem}

\proof Let $U$ be a connected component of
$V\in\Sigma^0_1(\mathbb{R}^n)$ and let $a$ be a rational point in
$U$. Then
 \begin{align*}
 U =\bigcup \big\{B(b,r) \mid\ & b\in \mathbb{Q}^n \wedge r \in \mathbb{Q}^+ \wedge
 \exists a_1 \cdots a_{m-1} \in \mathbb{Q}^n
                       \exists r_1\cdots r_{m-1} \in \mathbb{Q}^+\\
   &\big(\bigwedge_{i=1}^{m} C(a_i, r_i) \subseteq V
 \wedge \bigwedge_{i=1}^{m-1} (B(a_i, r_i) \cap B(a_{i+1},r_{i+1}) \neq \emptyset) \wedge a\in B(a_1,r_1)\big)\big\}
 \end{align*}
where $b=a_m$ and $r=r_m$. Since we can computably enumerate the basic
closed balls $C(b,r)\subseteq V$ \cite{kk07} with $b \in
\mathbb{Q}^n, r\in \mathbb{Q}^+$, $U$ is effectively open.
 \qed

\subsection{Effectively open sets in Euclidean spaces}

Now we prove the main result of this paper:

\begin{thm}\label{euclide}
For any $n\geq1$,
$FO(\Sigma^0_1(\mathbb{R}^n))\equiv_mFO(\mathbb{N})$.
\end{thm}

Since the upper bound holds by  Proposition \ref{loc}, we only
have to prove the lower bound
$FO(\mathbb{N})\leq_mFO(\Sigma^0_1(\mathbb{R}^n))$. Since
$\Sigma^0_1(\mathbb{R}^n)$ is a distributive lattice, we can use in
the definitions not only the symbol of inclusion but also the
symbols of Boolean operations and the constants
$\emptyset,\mathbb{R}^n$. 

Recall that a point  $x=(x_1,\ldots,x_n)\in \mathbb{R}^n$ is computable iff any real $x_i,i=1,\dots,n$, is computable, i.e. $x_i=\lim_np_n=\lim_nq_n$ for some computable sequences $\{p_n\},\{q_n\}$ of rationals with $p_0<p_1<\cdots x_i\cdots<q_1<q_0$. Thus, $x$
is computable iff $\mathbb{R}^n\setminus\{x\}$ is effectively open,
hence we can use (to simplify notation) in our defining formulas the
computable points (as the complements of effectively open sets maximal w.r.t. inclusion
among the effectively open sets strictly below $\mathbb{R}^n$). 

More precisely, we can use
a variable $x$ to range over the computable points,  the atomic
formulas $x\in U$ (where $U$ range as usual through
$\Sigma^0_1(\mathbb{R}^n)$), and we can quantify over $x$. Our
proof for $n\geq2$  is closely related to that in \cite{gr51} while
the proof for $n=1$ is based on quite different ideas, so we
consider the two cases separately. To keep our notation as close as possible to that of
\cite{gr51}, we let our formulas  use, along with the ``usual''
variables $U,V,W,G$ ranging over $\Sigma^0_1(\mathbb{R}^n)$,
the variables $a,x,y$ ranging over the computable points of
$\mathbb{R}^n$, and variables $A,B,C,D$ ranging over
$\Pi^0_1(\mathbb{R}^n)$. 

{\em Proof for $n\geq2$}. In this case we find an interpretation of
$\mathbb{N}$ in $\Sigma^0_1(\mathbb{R}^n)$ without
	parameters, which is sufficient by Lemma \ref{inter}. Let $\widetilde{\mathit{Cof}}$ be the unary relation on
$\Sigma^0_1(\mathbb{R}^n)$ which is true precisely on the
complements in $\mathbb{R}^n$ of finite sets of computable points (equivalently, we will use the relation $\widetilde{Fin}$ on $\Pi^0_1(\mathbb{R}^n)$ which is
true precisely on the  finite sets of computable points). Furthermore, let $\mathit{CON}(U)$ and $\mathit{CMP}(U,V)$ denote formulas of signature $\{\subseteq\}$ such that for all $U,V\in\Sigma^0_1(\mathbb{R}^n)$ we have: $\Sigma^0_1(\mathbb{R}^n)\models \mathit{CON}(U)$ iff $U$ is connected, and $\Sigma^0_1(\mathbb{R}^n)\models \mathit{CMP}(U,V)$ iff $U$ is a connected component of $V$. With these $\widetilde{\mathit{Fin}},\mathit{CON},\mathit{CMP}$ at hand (we define them later in the proof), it is easy to interpret $\mathbb{N}$ in $\Sigma^0_1(\mathbb{R}^n)$ similarly to \cite{gr51} or to the proof of Theorem \ref{euclide2}.

Namely, first we show that Axiom A6 in \cite{gr51} holds also in the effective setting, i.e. for any finite disjoint
sets $A,B$ of computable points the following formula   is true in $\Sigma^0_1(\mathbb{R}^n)$:
 \begin{eqnarray*}
(A \cup B \subseteq U \wedge CON(U)) \to \exists V
   \exists W (A \subseteq  V \wedge B \subseteq W \\
   {}\wedge V \cap W =
\emptyset \wedge V \cup W \subseteq U \wedge
\mathit{CON}(V)\wedge \mathit{CON}(W)).
 \end{eqnarray*}
 Indeed by the proof of the topological version of Axiom 6, given $A,B$ we can find open sets $V,W$ with the specified properties; moreover, these sets are obtained as finite unions of rational open balls, so $V,W$ are effectively open and the effective version of Axiom 6 holds.

Let now $A\approx_GB$ be the ternary relation meaning that $A,B$ are
finite disjoint sets of computable points and $G$ is an effectively
open set such that
 $$A\cup B \subseteq G\wedge \forall H
(\mathit{CMP}(H,G) \to (|H\cap A| =1 \wedge | H\cap B| =1)).$$
 Then $A\approx_GB$ implies that $A,B$ are of the same cardinality, $|A|=|B|$.
Note that for any  finite sets $A,B$ of computable points we have:
$|A|=|B|$ iff $\Sigma^0_1(\mathbb{R}^n)\models Eq(A,B)$ where
$Eq(A,B)$ is the formula
  $$\widetilde{\mathit{Fin}}(A)\wedge \widetilde{\mathit{Fin}}(B)\wedge\exists C,D ,G(\widetilde{\mathit{Fin}} ( C) \wedge \widetilde{\mathit{Fin}}( D)
  \wedge C  \approx_G D \wedge C=A\setminus B \wedge D = B\setminus A).$$

As in \cite{gr51}, it is now straightforward to interpret the
structure $\mathbb{N}$ in the structure $(\Sigma^0_1(\mathbb{R}^n);
	\widetilde{\mathit{Cof}})$ without parameters by interpreting natural numbers as the
cardinalities of finite sets $A,B$  of computable points (i.e., by
taking the quotient-set of all such $A$ under the equivalence
relation $Eq$) and interpreting  $+,\times$ as follows:
$|A|+|B|=|C|$ iff
 $$\exists A^\prime,B^\prime(Eq(A^\prime,A)\wedge Eq(B^\prime,B)
 \wedge A^\prime\cap B^\prime=\emptyset\wedge Eq(A^\prime\cup B^\prime,C),$$
 and $|A|\times |B|=|C|$ iff
 $$\exists U \exists F (Eq(B,F) \wedge F,C \subseteq U \wedge \forall V
(\mathit{CMP}(V,U) \to (| V\cap F| =1 \wedge | V\cap C| = |  A| ))).
$$

It remains to define  $\widetilde{Fin}, \mathit{CON},\mathit{CMP}$ from the beginning of the proof.
We will use the formulas $\mathit{Con}(U)$ and $\mathit{Cmp}(U,V)$
from the proof of Theorem \ref{euclide2} which now have a different
meaning because the set variables range now over the effectively open sets rather than the open sets. If $\Sigma^0_1(\mathbb{R}^n)\models \mathit{Con}(U)$ we say that $U$
is {\em effectively connected}. Note that if two effectively open sets
$U,V$ are effectively connected and $U\cap V\not=\emptyset$ then
$U\cup V$ is also effectively connected. Note also that connected
effectively open sets are effectively connected. If
$\Sigma^0_1(\mathbb{R}^n)\models \mathit{Cmp}(U,V)$ we say that $U$ is an
{\em effective connected component} of $V$.

Let $\Phi(V)$ be the
formula $\forall x\in V\exists U\subseteq V(x\in U\wedge
\mathit{Cmp}(U,V))$ saying that any computable point of $V$
belongs to some effective component  of $V$. Since the computable points are dense in $\mathbb{R}^n$, $\Sigma^0_1(\mathbb{R}^n)\models \Phi(V)$ implies that $V$
is the union of its effective components.

Define the unary relation $\widetilde{Iso}$  on
$\Pi^0_1(\mathbb{R}^n)$ (an effective analogue of the relation $Iso$
from \cite{gr51}) as follows: $\widetilde{Iso}(A)$ iff
\begin{multline*}
 \exists
V(A\subseteq V\wedge\Phi(V)\wedge\forall U(\mathit{Cmp}(U,V)\to
\exists!x(x\in A\cap U)))
\\
 \text{ and }
  \forall U(U\cap
A\not=\emptyset\to\exists x(x\in U\cap A)).
\end{multline*}

Then  the relation $\widetilde{Iso}$ is
definable, $\widetilde{Iso}(A)$ implies that any point in $A$ is
computable and isolated (indeed, take a satisfying $V$, then any $a\in A$ is in a unique effective component $U$ of $V$; choose the unique computable $x\in U\cap A$; since the computable points are dense in $\mathbb{R}^n$, $A\cap U=\{x\}$, hence $x=a$ is isolated), any infinite $\Pi^0_1$-set satisfying
$\widetilde{Iso}$ is not bounded, and any finite set of computable
points satisfies $\widetilde{Iso}$.

The definition of $\widetilde{Fin}$ looks as follows: $\widetilde{Fin}(A)$ iff
 $$\widetilde{Iso}(A)\wedge\forall U(\forall a\in A\exists V(a\in V\subseteq\overline{U})\to \exists W(A\subseteq W\subseteq\overline{U}))$$
  (the second conjunct says that if the closure $Cl(U)$ of an effectively open set $U$ is disjoint with $A$ then $A$ can be separated from $U$ by an effectively open set $W$). From left to right this is easy: if $V_a\in\Sigma^0_1(\mathbb{R}^n)$ satisfies $a\in V_a\subseteq\overline{U}$ for each $a\in A$ and $A$ is finite then $W:=\bigcup\{V_a\mid a\in A\}$ is effectively open and separates $A$ from $U$.
  
For the other direction, it suffices to show that if  $\widetilde{Iso}(A)$ holds and $A$ is infinite then there is an effectively open set $U$ such that $\forall a\in A(a\not\in Cl(U))$ and $A\subseteq W\in\Sigma^0_1(\mathbb{R}^n)$ implies $W\cap U\neq\emptyset$. Since  $\widetilde{Iso}(A)$, $\overline{A}$ is effectively open, hence one can effectively enumerate the rational closed balls $C(b,q)$ contained in $\overline{A}$. Hence, there is a computable sequence $\{B(b_k,q_k)\}$ of all rational open balls such that $C(b_k,q_k)\subseteq\overline{A}$ and $q_k\leq1$; note that $\bigcup_kB(b_k,q_k)=\overline{A}=\bigcup_kC(b_k,q_k)$. For any $k<\omega$, let $F(k):=0$ if $d(b_k,0)\leq q_k$, and let $F(k)$ be the integer part of $d(b_k,0)-q_k$ otherwise; the function $F$ is computable. As is well known \cite{ro67}, there is a strictly increasing limitwise-monotone-computable function $t:\omega\to\{1,2,\ldots\}$ that dominates all computable functions  on $\omega$ (recall that $t$ is limitwise-monotone-computable if for some computable function $g:\omega\times\omega\to\omega$ we have $g(k,s)\leq g(k,s+1)$ and $t(k)=\lim_sg(k,s)$, and that $t$ dominates a function $f:\omega\to\omega$ if there is $k_0$ with $\forall k\geq k_0(f(k)<t(k))$).

We claim that the set $U:=\bigcup_kB(b_k,r_k)$, where $r_k:=q_k-\frac{1}{t(F(k))}$, has the desired properties (note that in the possible case $r_k\leq0$ the ``ball'' $B(b_k,r_k)$ is empty). Obviously, $U\subseteq\overline{A}$. Since $F$ is computable and $t$ is limitwise-monotone-computable, $U$ is effectively open. Let us check that $a\not\in Cl(U)$ for each $ a\in A$. Suppose the  contrary: there is a sequence $\{c_i\}$ in $U$ that converges to $a$, so in particular for some $i_0$ we have $\forall i\geq i_0(d(c_i,a)\leq1)$. Choose a sequence of naturals $\{k_i\}$ with $c_i\in B(b_{k_i},r_{k_i})$. Since $a\not\in B(b_{k_i},q_{k_i})$, we have $d(c_i,a)\geq\frac{1}{t(F(k_i))}$ for all $i$. Note that $\forall i(F(k_i)\leq m)$ for some $m$, because for all $i\geq i_0$ we have 
 $$F(k_i)\leq d(b_{k_i},0)\leq d(b_{k_i},c_i)+d(c_i,0)\leq q_{k_i}+d(c_i,0)\leq 1+d(c_i,a)+d(a,0)\leq1+1+d(a,0).$$
  For any $i$ we have $\frac{1}{t(F(k_i))}\geq \frac{1}{t(m)}$ and therefore $d(a,c_i)\geq\frac{1}{t(m)}>0$, contradicting the convergence of $\{c_i\}$ to $a$.

Finally, let $A\subseteq W\in\Sigma^0_1(\mathbb{R}^n)$; we have to show that $W\cap U\neq\emptyset$. Since $W\setminus A$ is effectively open and $A$ is infinite (hence unbounded), there is a computable subsequence $\{B(b_{k_i},q_{k_i})\}$ of $\{B(b_{k},q_{k})\}$ such that $B(b_{k_i},q_{k_i})\subseteq W$ for each $i$ and the sequence $\{F(k_i)\}$ is strictly increasing. Let $h:\omega\to\mathbb{Q}^+$ be a computable function such that $q_{k_i}=h(F(k_i))$ for each $i<\omega$. Since $t$ dominates all computable functions, for some $i_0$ we have $\forall m\geq i_0(\frac{1}{h(m)}<t(m))$. Since $\{F(k_i)\}$ is strictly increasing, $F(k_i)\geq i_0$ for some $i$, hence $\frac{1}{h(F(k_i))}<t(F(k_i))$ and therefore $q_{k_i}>\frac{1}{t(F(k_i))}$. Thus, the non-empty ball $B(b_{k_i},r_{k_i})$ is contained in $W\cap U$, so the latter set is non-empty. This completes the proof of definability of $\widetilde{Fin}$.

Let $\mathit{Bound}_\Sigma(U):=\forall A\subseteq U(\widetilde{Iso}(A)\to\widetilde{Fin}(A))$. Then for any  $U\in\Sigma^0_1(\mathbb{R}^n)$ we have: $\Sigma^0_1(\mathbb{R}^n)\models Bound_\Sigma(U)$ iff $U$ is bounded. Indeed, from left to right this follows from the properties of $\widetilde{Iso}$. Conversely, let $U$ be non-bounded.  Then there is a computable sequence $\{a_k\}$ of rational points in $U$ such that $d(a_{k+1},0)\geq d(a_k,0)+1$ for each $k<\omega$. Then  $A:=\{a_k\mid k<\omega\}$ is an infinite $\Pi^0_1$-subset of $U$ and $\widetilde{Iso} (A)$ holds, as desired.

Let $\mathit{Bound}_\Pi(A):=\exists U\supseteq A(Bound_\Sigma(U))$. Then clearly for any  $A\in\Pi^0_1(\mathbb{R}^n)$ we have: $\Sigma^0_1(\mathbb{R}^n)\models Bound_\Pi(A)$ iff $A$ is bounded.

Let
 $$\mathit{BCON}(A):=\mathit{Bound}_\Pi(A)\wedge\neg\exists U,V(A\subseteq U\cup V\wedge U\cap V=\emptyset\wedge A\cap U\neq\emptyset\wedge A\cap V\neq\emptyset)).$$
  Then  for any  $A\in\Pi^0_1(\mathbb{R}^n)$ we have: $\Sigma^0_1(\mathbb{R}^n)\models BCON(A)$ iff $A$ is bounded and connected. Indeed, from right to left this is obvious. Conversely, it suffices to show that if $A$ bounded and non-connected then
   $$\Sigma^0_1(\mathbb{R}^n)\models \exists U,V(A\subseteq U\cup V\wedge U\cap V=\emptyset\wedge A\cap U\neq\emptyset\wedge A\cap V\neq\emptyset)).$$
    Since $A$ in non-connected, for some open sets $U',V'$ we have
     $$A\subseteq U'\cup V', U'\cap V'=\emptyset, A\cap U'\neq\emptyset, A\cap V'\neq\emptyset.$$
      Representing $U'$ and $V'$ as the unions of some families $\mathcal{F}'$ and $\mathcal{G}'$ of rational open balls, we obtain an open cover $\mathcal{F}'\cup\mathcal{G}'$ of $A$. Since $A$ is compact, there are finite subfamilies $\mathcal{F}$ and $\mathcal{G}$ of resp. $\mathcal{F}'$ and $\mathcal{G}'$ such that $\mathcal{F}\cup\mathcal{G}$ covers $A$. Then the effectively open sets $U:=\bigcup\mathcal{F}$ and $V:=\bigcup\mathcal{G}$ have the desired properties.

Let
 $$\mathit{CON}(U):=\forall x,y\in U\exists B(x,y\in B\subseteq U\wedge \mathit{BCON}(B)).$$
  Then for any  $U\in\Sigma^0_1(\mathbb{R}^n)$ we have: $\Sigma^0_1(\mathbb{R}^n)\models \mathit{CON}(U)$ iff $U$ is connected. Indeed, if $U$ is connected and  $x,y$ are computable points in $U$, let $B$ be a polygonal path with rational inner points that connects $x$ and $y$ inside $U$; then $x,y\in B\subseteq U$ and $\mathit{BCON}(B)$, as desired. Conversely, let $U$ be non-connected. Choose computable points $x,y$  from distinct connected components of $U$. Then clearly there is no connected set $B$ with $x,y\in B\subseteq U$. 

With the formula $\mathit{CON}$ at hand, it is straightforward to write down the formula $\mathit{CMP}$, completing the proof of the theorem for the case $n>1$.

In fact, the arguments above (except those concerning Axiom 6) work for any $n\geq1$, so we have the following corollary which is interesting in its own right.

\begin{cor}\label{connect}
There are formulas  $\mathit{CON}(U)$ and $\mathit{CMP}(U,V)$ of signature $\{\subseteq\}$ such that for all $n\geq1$ and $U,V\in\Sigma^0_1(\mathbb{R}^n)$ we have: $\Sigma^0_1(\mathbb{R}^n)\models \mathit{CON}(U)$ iff $U$ is connected, and $\Sigma^0_1(\mathbb{R}^n)\models \mathit{CMP}(U,V)$ iff $U$ is a connected component of $V$.
\end{cor}

\noindent {\em Proof of Theorem \ref{euclide} for $n=1$}.
We use the formulas $\mathit{CON}$ and $\mathit{CMP}$ from Corollary \ref{connect}.  Let
$\xi(U,V)$ be the formula
 $$
U\not=\emptyset\wedge V\not=\emptyset\wedge U\cap V=\emptyset\wedge
\forall U^\prime(U^\prime\cap V=\emptyset\rightarrow
U^\prime\subseteq U)
 \wedge\\ \forall V^\prime(V^\prime\cap U=\emptyset\rightarrow V^\prime\subseteq V)
 $$
 saying in $\Sigma^0_1(\mathbb{R})$ that $U,V$ are disjoint non-empty effectively open sets such that
$U=Int(\mathbb{R}\setminus V)$ (where $Int$ is the interior
operator) and $V=Int(\mathbb{R}\setminus U)$.

Note that between any $U$-components $U_1<U_2$ (where $U_1<U_2$ means $\forall x\in U_0\forall y\in
U_1(x<y)$)) there
is a $V$-component. Suppose the contrary, then the interval
$W:=(inf(U_1),sup(U_2))$ is disjoint with $V$, hence $W\subseteq U$.
This is a contradiction because $sup(U_1)\in W\setminus U$. Symmetrically,  between any two $V$-components there is a $U$-component.

Let
$\{q_0,q_1,\cdots\}$ be a computable enumeration of the set
$U\cap\mathbb{Q}$ without repetitions. Define the following equivalence
relation $\sim$ on $\omega$: $m\sim n$ iff $q_m,q_n$ are in the
same $U$-component. We claim that this relation is c.e. Indeed,
since $m\sim n$ is equivalent to the disjunction of $q_m\leq
q_n\wedge[q_m,q_n]\subseteq U$ and $q_n\leq
q_m\wedge[q_n,q_m]\subseteq U$, so it suffices to check that the
relation $q_m\leq q_n\wedge[q_m,q_n]\subseteq U$ is c.e. We have
$U=\bigcup_iB_i$ for a computable sequence $\{B_i\}$ of basic open
balls (i.e., intervals with rational endpoints). Since closed
intervals are compact, the relation $q_m\leq
q_n\wedge[q_m,q_n]\subseteq U$ is equivalent to
 $$\exists l\exists i_0,\ldots,i_l(q_m\in B_{i_0}\wedge q_n\in B_{i_l}
 \wedge B_{i_0}\cap  B_{i_1}\not=\emptyset \wedge  \cdots \wedge B_{i_{l-1}}\cap  B_{i_l}\not=\emptyset).$$
 Alternatively, the last assertion follows from the results in
 \cite{kk07}.
Therefore, $\sim$ is c.e. It is also co-c.e. because, by the
previous paragraph, $m\not\sim n$ is equivalent to the disjunction
of the predicates $q_m<q_n\wedge\exists r\in V\cap\mathbb{Q}
(q_m<r<q_n)$ and $q_n<q_m\wedge\exists r\in V\cap\mathbb{Q}
(q_n<r<q_m)$.

Let now
 $$\mathit{Cmp}^\ast(U,V):=U\subseteq V\wedge\forall U^\prime(CMP(U^\prime,V)\rightarrow
 (U^\prime\subseteq U\vee\ U\cap U^\prime=\emptyset)).$$
 Then
$\Sigma^0_1(\mathbb{R})\models Cmp^\ast(U,V)$ iff $U$ is an effectively open union
of some connected components of $V$. 
Let
 $$\mathit{ICmp}(V):=\exists U(\mathit{Cmp}^\ast(U,V)\wedge\neg\exists W(W=V\setminus U\wedge
 \mathit{Cmp}^\ast(W,V))).$$
 Then $\Sigma^0_1(\mathbb{R})\models \mathit{ICmp}(V)$ implies that  $V$ has
infinitely many connected components. Suppose the contrary, then
$V=V_0\cup\cdots\cup V_n$ for some $n\geq 0$ and pairwise disjoint
components $V_0,\ldots,V_n \in\Sigma^0_1(\mathbb{R})$. Then any $U$
with $\Sigma^0_1(\mathbb{R})\models \mathit{Cmp}^\ast(U,V)$ is a union of
some of $V_0,\ldots,V_n$, hence $V\setminus U$ is the finite union
of the remaining $V_i$, a contradiction.

Finally, let
 $\theta(U,V):=
 \xi(U,V)\wedge \mathit{ICmp}(U).
 $
 Then $\Sigma^0_1(\mathbb{R})\models \theta(U,V)$ iff both $U,V$ have infinitely
many  connected components which are computable intervals,
$U=Int(\mathbb{R}\setminus V)$ and $V=Int(\mathbb{R}\setminus U)$.
For any such $U,V$,  the relation $\sim$ is computable with infinitely many
equivalence classes,  hence the lattice $\mathcal{E}$ is isomorphic
to the lattice $\mathcal{F}$ of c.e. sets closed under $\sim$. The
lattice $\mathcal{F}$ is in turn isomorphic to the lattice
$(\mathcal{G};\subseteq)$ where
$\mathcal{G}:=\{G\mid\Sigma^0_1(\mathbb{R})\models \mathit{Cmp}^\ast(G,U)\}$
consists of effectively open subsets of $U$ closed under components,
for each $U$ as above.
This shows that $\mathcal{E}$ is $c$-definable in $\Sigma^0_1(\mathbb{R})$ with a non-empty definable set of parameters. By Lemma \ref{inter}, $FO(\mathcal{E})\leq_mFO(\Sigma^0_1(\mathbb{R}))$.
 \qed

\subsection{Effectively open sets in domains}

 For some natural domains we have the following interpretability result:

\begin{thm}\label{scott}
 For any $X\in\{P\omega,k^{\leq\omega},
 k^\omega_\bot\mid 2\leq k\leq\omega\}$,  $\mathcal{E}$ is
 interpretable without parameters in  $\Sigma^0_1(X)$.
\end{thm}

\proof  We give the proof only for the most popular space $P\omega$ (which is homeomorphic to $1^\omega_\bot$) but
similar arguments  work for the other spaces as well. First we check
(just as in the proof of Theorem \ref{scott2}) that the set of
finitary elements $\orduparrow{F}$, for all finite $F\subseteq\omega$,
is definable in $\Sigma^0_1(P\omega)$.

Next  we check (again as in the proof of Theorem \ref{scott2}) that
the set $\mathcal{V}_n:=\{\orduparrow{F}:|F|=n\}$ is definable in
$\Sigma^0_1(P\omega)$ for each $n<\omega$ via the formula
$\phi_n(V)$.

Now let $U_n:=\{S\subseteq\omega:n\leq|S|\}$, so $U_0=P\omega$ and
$U_n=\bigcup \mathcal{V}_n$. Then the singleton set $\{U_{n+1}\}$ is
defined by the formula
 $$\psi_{n+1}(U):=\forall V(\phi_{n+1}(V)\rightarrow V\subseteq U)\wedge
 \neg\exists V(\phi_n(V)\wedge V\subseteq u).$$
 By the proof of Theorem \ref{metric}, the lattice $\mathcal{E}$ is
isomorphic to the sublattice $\{S\in\Sigma^0_1(P\omega)\mid
U_1\subseteq S\subseteq U_2\}$ of $\Sigma^0_1(P\omega)$, and is thus
definable without parameters.
 \qed

\begin{cor}\label{scott1}
For any $X\in\{P\omega,k^{\leq\omega},
 k^\omega_\bot\mid 2\leq k\leq\omega\}$, $FO(\Sigma^0_1(X))\equiv_mFO(\mathbb{N})$.
\end{cor}

\proof The  upper bound holds by Proposition \ref{loc}, the
lower bound by the previous theorem (alternatively, by the proof of
Theorem \ref{scott2}).
 \qed

\section{Conclusion}\label{con1}

The  results of this paper show that a satisfactory understanding of definability
in the lattices of open and of effectively open sets is connected
with intricate relationships between topological and algorithmic
properties of the corresponding effective spaces. For this reason we
believe that this research direction is interesting and deserves  further developments.
Many questions remain open. In particular, we are still far
from understanding the border between  decidability and undecidability of the
theories of  lattices of all open sets: in particular, we 
currently have no general sufficient condition giving undecidability of
this theory (similar to Theorem \ref{metric} for the effectively
open sets).

The methods  in Section \ref{efopen} are rather different when we distinguish between 
metric spaces and  domains. It would be useful to develop  unified
methods applicable e.g. to many quasi-metric spaces \cite{br13}.

The methods of  this paper work mainly for  second countable locally
compact spaces. It would be nice to investigate our questions for
 second countable  non-locally compact spaces like the Baire space or for
the space $\mathbb{R}^\omega$. The situation with non- second countable 
spaces is even less clear.

We guess that Theorem \ref{scott} and Corollary \ref{scott1} hold
also for the spaces of continuous partial functionals of finite
types over $\omega$ but the given proofs should be modified
considerably.

A natural generalization of the topic of this paper is the study of
first order theories of lattices of other levels of the standard
hierarchies, in particular of higher levels of Borel hierarchy
$\mathbf{\Sigma}^0_n$ and of the effective Borel hierarchy
$\Sigma^0_n$. As shown in \cite{ra69}, the first order theories of
the lattices of $\mathbf{\Sigma}^0_2$-sets in $\mathcal{C}$ and in
$\mathbb{R}$ are decidable.

\section*{Acknowledgement}
The first author was supported by RFBR project 17-01-00247-a. 
The second author was supported by the project OpenLab of Kazan (Volga region) Federal University. Both authors received funding from
the People Programme (Marie Curie Actions) of the European Union's Seventh Framework
Programme FP7/2007-2013 under the REA grant agreement no.\
PIRSES-GA-2011-294926-COMPUTAL.

The authors  thank Andr\'e Nies for a discussion of the
lower bound problem for the theory $FO(\Sigma^0_1(\mathcal{N}))$ and the three anonymous referees for many helpful suggestions.

\end{document}